\documentclass[10pt,twoside,reqno]{amsart}
\usepackage{amssymb}
\usepackage{multirow}
\calclayout

\numberwithin{equation}{section}
\makeatletter

\renewcommand{\@secnumfont}{\bfseries}

\renewcommand{\section}{\@startsection{section}{1}%
  {0mm}{.7\linespacing\@plus\linespacing}{.5\linespacing}
  {\normalfont\bfseries\centering}}

\newcommand{\bibsection}{\@startsection{section}{1}%
  {0mm}{.7\linespacing\@plus\linespacing}{.5\linespacing}
  {\normalfont\scshape\centering}}

\renewcommand{\@biblabel}[1]{#1.}

\newtheorem{thm}{\bf Theorem}[section]
\newtheorem{lem}[thm]{\bf Lemma}
\newtheorem{cor}[thm]{\bf Corollary}

\begin{document}

\vspace{1.3cm}

\title {Some identities involving special numbers and moments of random variables}

\author{Taekyun Kim}
\address{Department of Mathematics, College of Science, Tianjin Polytechnic University, Tianjin City, 300387, China
\\ Department of Mathematics, Kwangwoon University, Seoul 139-701, Republic
	of Korea}
\email{tkkim@kw.ac.kr}

\author{Yonghong Yao}
\address{Institute of Fundamental Sciences, University of Electronic Science and Technology, Chengdu 611731,
 China \\Department of Mathematics, College of Science, Tianjin Polytechnic University, Tianjin City, 300387, China}
\email{yaoyonghong@aliyun.com}

\author{Dae San Kim}
\address{Department of Mathematics, Sogang University, Seoul 121-742, Republic of Korea}
\email{dskim@sogang.ac.kr}

\author{Hyuck-In Kwon}
\address{Department of Mathematics, Kwangwoon University, Seoul 139-701, Republic of Korea}
\email{sura@kw.ac.kr}

\subjclass[2010]{05A19; 11B83}
\keywords{Random variable, Moment, Stirling number, Degenerate Stirling number, Derangement number}
\begin{abstract}
In this paper, we derive some identities involving special numbers and moments of random variables by using the generating functions of the moments of certain random variables. Here the related special numbers are Stirling numbers of the first and second kinds, degenerate Stirling numbers of the first and second kinds, derangement numbers, higher-order Bernoulli numbers and Bernoulli numbers of the second kind.
\end{abstract}
\maketitle

\markboth{\centerline{\scriptsize Some identities involving special numbers and moments of random variables }}
{\centerline{\scriptsize T. Kim, Y. Yao,  D. S. Kim, H.-I. Kwon}}

\bigskip
\medskip

\section{Introduction}
As is well known, the Bernoulli polynomials of order $r$ are defined by the generating function
\begin{equation}\begin{split}\label{01}
\left( \frac{t}{e^t-1} \right)^r e^{xt} = \sum_{n=0}^\infty  B_n^{(r)}(x)
 \frac{t^n}{n!},\quad (\textnormal{see} \,\, [9,10]).
\end{split}\end{equation}
When $x=0$, $B_n^{(r)}=B_n^{(r)}(0)$ are called the Bernoulli numbers of order $r$. In particular, for $r=1$, $B_n(x)=B_n^{(1)}(x)$ are the Bernoulli polynomials and $B_n=B_n^{(1)}$ are the Bernoulli numbers.

For $\lambda \in \mathbb{R}$, L. Carlitz considered the degenerate Bernoulli polynomials of order $r$ which are given by
\begin{equation}\begin{split}\label{02}
\left( \frac{t}{(1+\lambda t)^{\frac{1}{\lambda }}-1} \right)^r (1+\lambda t)^{\frac{x}{\lambda }} = \sum_{n=0}^\infty         \beta_{n,\lambda }^{(r)}(x)         \frac{t^n}{n!},\quad (\textnormal{see} \,\, [1,2]).
\end{split}\end{equation}
When $x=0$, $\beta_{n,\lambda }^{(r)}=\beta_{n,\lambda }^{(r)}(0)$ are called the degenerate Bernoulli number of order $r$. In particular, for $r=1$, $\beta_{n,\lambda }^{(1)}(x)=\beta_{n,\lambda }(x)$ are the degenerate Bernoulli polynomials.

For $n \in \mathbb{N}$, the falling factorial sequence is defined as
\begin{equation*}\begin{split}
(x)_0=1 ,\,\,(x)_n = x(x-1)\cdots(x-n+1), (n\geq1).
\end{split}\end{equation*}

The Stirling numbers of the first kind are defined as
\begin{equation*}\begin{split}
S_1(0,0)=1,\,\,(x)_n = \sum_{l=0}^n S_1(n,l) x^l,\,\,(n \geq 1),\quad (\textnormal{see} \,\, [5,9]).
\end{split}\end{equation*}

The Stirling numbers of the second kind are given by
\begin{equation*}\begin{split}
x^n = \sum_{l=0}^n S_2(n,l) (x)_l, \,\,(n \geq 0 ).
\end{split}\end{equation*}

In [6], the degenerate Stirling numbers of the first kind are introduced by the generating function
\begin{equation}\begin{split}\label{03}
\frac{1}{k!} \Big( \log(1+\lambda t)^{\frac{1}{\lambda }} \Big)^k = \sum_{n=k}^\infty S_{1,\lambda }(n,k) \frac{t^n}{n!},\,\,(\lambda  \in \mathbb{R}).
\end{split}\end{equation}

From \eqref{03}, we have
\begin{equation}\begin{split}\label{04}
(x)_{n,\lambda } = \sum_{l=0}^n S_{1,\lambda }(n,l) x^l,\,\,(n \geq 0),\quad (\textnormal{see} \,\, [5,6]).
\end{split}\end{equation}
where $(x)_{n,\lambda }=x(x-\lambda )\cdots(x-(n-1)\lambda ), (n \geq 1)$, $(x)_{0,\lambda }=1$.

By \eqref{04}, we get
\begin{equation}\begin{split}\label{05}
S_{1,\lambda }(n+1,k) = S_{1,\lambda }(n,k-1)-n\lambda S_{1,\lambda }(n,k) ,\,\,(1 \leq k \leq n),\quad (\textnormal{see} \,\, [6]).
\end{split}\end{equation}

In [5], the degenerate Stirling numbers of the second kind are defined by
\begin{equation}\begin{split}\label{06}
\frac{1}{k!} \Big( (1+\lambda t)^{\frac{1}{\lambda }}-1 \Big)^k = \sum_{n=k}^\infty S_{2,\lambda }(n,k) \frac{t^n}{n!},\,\,(k \geq 0).
\end{split}\end{equation}

From \eqref{06}, we have
\begin{equation}\begin{split}\label{07}
\sum_{m=0}^n \lambda ^{n-m} S_1(n,m) \frac{1}{k!}\Delta^k0^m = \begin{cases}
S_{2,\lambda }(n,k),&\text{if}\,\, n \geq k,\\
0,&\text{otherwise},
\end{cases}
\end{split}\end{equation}

and
\begin{equation}\begin{split}\label{08}
S_{2,\lambda }(n+1,k) = kS_{2,\lambda }(n,k) + S_{2,\lambda }(n,k-1) - n\lambda S_{2,\lambda }(n,k),
\end{split}\end{equation}
where $\Delta f(x)=f(x+1)-f(x)$, and  $1 \leq k \leq n$. Note that, letting $\lambda \rightarrow 0$ gives us
\begin{equation*}
S_2(n+1,k)=k S_2(n,k) + S_2(n,k-1).
\end{equation*}

The Bernoulli numbers of the second kind are defined as
\begin{equation}\begin{split}\label{09}
\frac{t}{\log(1+t)} = \sum_{n=0}^\infty b_n \frac{t^n}{n!},\quad (\textnormal{see} \,\, [3,8,9]).
\end{split}\end{equation}

It is known that
\begin{equation}\begin{split}\label{10}
\left( \frac{t}{\log(1+t)} \right)^k (1+t)^{x-1} = \sum_{n=0}^\infty B_n^{(n-k+1)}(x) \frac{t^n}{n!},\quad (\textnormal{see} \,\, [9]).
\end{split}\end{equation}

From \eqref{10}, we note that
\begin{equation*}\begin{split}
b_n = B_n^{(n)}(1),\,\,(n \geq 0).
\end{split}\end{equation*}

A random variable $X$ is a real-valued function defined on a sample space. We say that $X$ is a continuous random variable if there exists a nonnegative function $f(x)$, defined for all $x \in (-\infty, \infty)$, having the property that for any set $B$ of real numbers
\begin{equation}\begin{split}\label{11}
P\{ X \in B \} = \int_B f(x) dx,\quad (\textnormal{see} \,\, [9]).
\end{split}\end{equation}

The function $f(x)$ is called the probability density function of the random variable $X$.

Let $X$ be a uniform random variable on the interval $(\alpha, \beta)$. Then the probability density function of $X$ is given by
\begin{equation}\begin{split}\label{12}
f(x) = \begin{cases}
\frac{1}{\beta-\alpha},& \text{if}\,\,\alpha<x<\beta,\\
0,&\text{otherwise}.
\end{cases}
\end{split}\end{equation}
A continuous random variable whose density function is given by
\begin{equation}\begin{split}\label{12}
f(x) = \begin{cases}
\frac{\lambda e^{-\lambda x}(\lambda x)^{\alpha-1}}{\Gamma(\alpha)},& \text{if}\,\,x \geq 0,\\
0,&\text{if}\,\, x <0,
\end{cases}
\end{split}\end{equation}
for some $\lambda >0$, $\alpha>0$ is said to be the gamma random variable with parameter $\alpha, \lambda$ (see [9]).

If $X$ is a continuous random variable having a probability density function $f(x)$, then the expectation of $X$ is defined by
\begin{equation}\begin{split}\label{14}
E[X]=\int_{-\infty}^\infty x f(x) dx.
\end{split}\end{equation}

Let $X$ be continuous random variable with the probability density function $f(x)$. For any real-valued function $g$, we have
\begin{equation}\begin{split}\label{15}
E[g(X)] =\int_{-\infty}^\infty g(x) f(x) dx,\quad (\textnormal{see} \,\, [9]).
\end{split}\end{equation}

The expected value of a random variable $X$, $E[X]$, is also referred to as the mean or the first moment of $X$. The quantity $E[X^n]$, $n \geq 1$, is said to be the $n$-th moment of $X$. That is,
\begin{equation}\begin{split}\label{16}
E[X^n] = \int_{-\infty}^\infty x^n f(x) dx,\,\,(n \geq 1).
\end{split}\end{equation}

Another quantity of interest is the variance of random variable $X$ which is defined by
\begin{equation}\begin{split}\label{17}
Var(X) = E[ (X-E[X])^2]= E[X^2]- (E[X])^2.
\end{split}\end{equation}

We say that $X$ and $Y$ are jointly continuous if there exists a function $f(x,y)$, defined for all $x$ and $y$, having the property that for all sets $A$ and $B$ of real numbers
\begin{equation}\begin{split}\label{18}
P\{X \in A, \,\,Y \in B\} = \int_B \int_A f(x,y) dxdy.
\end{split}\end{equation}

The function $f(x,y)$ is called the joint probability density function of $X$ and $Y$.

Let $f_X(x)$ be the probability density function of $X$. Then we have
\begin{equation*}\begin{split}
f_X(x) = \int_{-\infty}^\infty f(x,y) dy, \,\,f_Y(y) = \int_{-\infty}^\infty f(x,y) dx,\quad (\textnormal{see} \,\, [9]).
\end{split}\end{equation*}

The random variables $X$ and $Y$ are independent if
\begin{equation}\begin{split}\label{19}
P\{ X \leq a, \,\, Y \leq b \} = P\{ X \leq a \} \cdot P \{ Y \leq b \}.
\end{split}\end{equation}

Let $X$ and $Y$ be independent random variables. For any real-valued functions $h$ and $g$, we have
\begin{equation}\begin{split}\label{20}
E[g(X)h(Y)] = E[g(X)] E[h(Y)],\quad (\textnormal{see} \,\, [9]).
\end{split}\end{equation}
and
\begin{equation}\begin{split}\label{21}
f(x,y) = f_X(x)f_Y(y).
\end{split}\end{equation}

Let $X,Y$ have a joint probability density function $f(x,y)$. Then the conditional probability density function of $X$, given that $Y=y$, is defined for all values of $y$ such that $f_Y(y) >0$, by
\begin{equation}\begin{split}\label{21}
f_{X|Y} (x|y) = \frac{f(x,y)}{f_Y(y)}.
\end{split}\end{equation}

The conditional expectation of $X$, given that $Y=y$, is defined for all values of $y$ such that $f_Y(y) >0$, by
\begin{equation}\begin{split}\label{22}
E[X|Y=y] = \int_{-\infty}^\infty x f_{X|Y} (x|y)dx,\quad (\textnormal{see} \,\, [9]).
\end{split}\end{equation}

Let $E[X|Y]$ be the function of random variable $Y$ whose value at $Y=y$ is $E[X|Y=y]$.

Then, we note that
\begin{equation}\begin{split}\label{23}
E[X] = E[E[X|Y]] = \int_{-\infty}^\infty E[X|Y=y] f_Y(y) dy.
\end{split}\end{equation}

In this paper, we derive some identities involving special numbers and moments of random variables by using the generating functions of the moments of certain random variables. Here the related special numbers are Stirling numbers of the first and second kinds, degenerate Stirling numbers of the first and second kinds, derangement numbers, higher-order Bernoulli numbers and Bernoulli numbers of the second kind.

\section{Explicit formulas arising from probabilistic representations}

Let $U_1, U_2 , \cdots, U_k$ be uniformly independent random variables on (0,1). Then we have

\begin{equation}\begin{split}\label{24}
&E[e^{(U_1+U_2+\cdots+U_k)(e^t-1)}] = E[e^{U_1(e^t-1)}] \cdots E[e^{U_k(e^t-1)}]\\
&=\int_0^1 e^{u_1(e^t-1)}f(u_1)du_1 \times \int_0^1 e^{u_2(e^t-1)}f(u_2)du_2 \times \cdots \times \int_0^1 e^{u_k(e^t-1)}f(u_k)du_k\\
&=\left( \frac{t}{e^t-1} \right)^k \frac{k!}{t^k} \frac{1}{k!} \left( e^{(e^t-1)}-1\right)^k =\left( \frac{t}{e^t-1} \right)^k \frac{k!}{t^k}  \sum_{l=k}^\infty S_2(l,k) \frac{1}{l!} (e^t-1)^l \\
&= \left( \sum_{j=0}^\infty B_j^{(k)} \frac{t^j}{j!} \right) \times \frac{k!}{t^k} \left( \sum_{m=k}^\infty \sum_{l=k}^m S_2(l,k) S_2(m,l) \frac{t^m}{m!} \right)\\
&= \left( \sum_{j=0}^\infty B_j^{(k)} \frac{t^j}{j!} \right) \left( \sum_{m=0}^\infty \left( \sum_{l=k}^{m+k} S_2(l,k) S_2(m+k,l) \frac{k! m!}{(m+k)!} \right) \frac{t^m}{m!} \right) \\
&=\sum_{n=0}^\infty \left( \sum_{ m=0}^n \sum_{l=k}^{m+k} \frac{{n \choose m}}{{m+k \choose m}} S_2(l,k) S_2(m+k,l) B_{n-m}^{(k)} \right) \frac{t^n}{n!}.
\end{split}\end{equation}

On the other hand,
\begin{equation}\begin{split}\label{25}
&E[e^{(U_1+U_2+\cdots+U_k)(e^t-1)}] = \sum_{m=0}^\infty E[(U_1+\cdots+U_k)^m] \frac{(e^t-1)^m}{m!}\\
&= \sum_{m=0}^\infty E[(U_1+\cdots+U_k)^m] \sum_{n=m}^\infty S_2(n,m) \frac{t^n}{n! }\\
&= \sum_{n=0}^\infty \left( \sum_{m=0}^n S_2(n,m)  E[(U_1+\cdots+U_k)^m]  \right) \frac{t^n}{n!}.
\end{split}\end{equation}

Therefore, by \eqref{24} and \eqref{25}, we obtain the following theorem.

\begin{thm}
Let $U_1, U_2, \cdots, U_k$ be uniformly independent random variables on (0,1). For $k \geq 0$, we have
\begin{equation*}\begin{split}
 \sum_{ m=0}^n \sum_{l=k}^{m+k} \frac{{n \choose m}}{{m+k \choose m}} S_2(l,k) S_2(m+k,l) B_{n-m}^{(k)} = \sum_{m=0}^n S_2(n,m)  E[(U_1+\cdots+U_k)^m] .
\end{split}\end{equation*}
\end{thm}

\begin{cor}
For $n,k \geq 0$, we have
\begin{equation*}\begin{split}
 &\sum_{ m=0}^n \sum_{l=k}^{m+k} \frac{{n \choose m}}{{m+k \choose m}} S_2(l,k) S_2(m+k,l) B_{n-m}^{(k)}\\
 & = \sum_{m=0}^n \sum_{l_1+\cdots+l_k=m+k, l_i \geq 1} {m \choose l_1,l_2, \cdots,l_k} S_2(n,m).
\end{split}\end{equation*}
\end{cor}

Let $U$ be a uniform random variable on (0,1). Then we have
\begin{equation}\begin{split}\label{26}
E[(1+t)^U] &= \int_0^1 (1+t)^u p(u) du = \int_0^1 (1+t)^x dx = \frac{t}{\log(1+t)}\\
&=\sum_{n=0}^\infty        b_n          \frac{t^n}{n!}.
\end{split}\end{equation}

On the other hand,
\begin{equation}\begin{split}\label{27}
E[(1+t)^U] &= E[e^{U \log(1+t)} ] = \sum_{k=0}^\infty E[U^k] \frac{1}{k!} \log^k(1+t) \\
&= \sum_{k=0}^\infty E[U^k] \sum_{n=k}^\infty S_1(n,k) \frac{t^n}{n!} = \sum_{n=0}^\infty \left( \sum_{k=0}^n E[U^k] S_1(n,k) \right) \frac{t^n}{n!}.
\end{split}\end{equation}

Thus, by \eqref{26} and \eqref{27}, we easily get
\begin{equation*}\begin{split}
b_n = \sum_{k=0}^n E[U^k] S_1(n,k),\,\,(n \geq k \geq 0).
\end{split}\end{equation*}

\begin{lem}
Let $U$ be a uniform random variable on (0,1). For $n,k \geq 0$, we have
\begin{equation*}\begin{split}
b_n = \sum_{k=0}^n E[U^k] S_1(n,k).
\end{split}\end{equation*}
\end{lem}

Let $X$ be a gamma random variable with parameters $\alpha=u(>0)$ and $\lambda=1$ and let $U$ be a uniform random variable on (0,1). Assume that $X$ and $U$ are independent. Then we have
\begin{equation}\begin{split}\label{28}
E[e^{Xt}] &= E[E[e^{Xt}|U]] = \int_0^1 E[e^{Xt}|U=u] f(u) du \\
&= \int_0^1 \int_0^\infty e^{xt} f_{X|U} (x|u) dxdu = \int_0^1 \int_0^\infty e^{xt} \frac{f_X(x) f_U(u)}{f_U(u)} dx du\\
&=  \int_0^1  \frac{1}{\Gamma(u)}\int_0^\infty e^{xt} e^{-x} x^{u-1} dx du =  \int_0^1 \left( \frac{1}{1-t} \right)^u du \\
&= \frac{-t}{(1-t)\log(1-t)} ,\,\,( 0<t<1).
\end{split}\end{equation}

From \eqref{28}, we note that
\begin{equation}\begin{split}\label{29}
(1-t)E[e^{Xt}] = \frac{-t}{\log(1-t)} = \sum_{n=0}^\infty b_n (-1)^n \frac{t^n}{n!}.
\end{split}\end{equation}

We observe that
\begin{equation}\begin{split}\label{30}
(1-t)E[e^{Xt}] &=  \sum_{n=0}^\infty E[X^n] \frac{t^n}{n!} - \sum_{n=0}^\infty                  E[X^n]\frac{t^{n+1}}{n!} \\
&= \sum_{n=0}^\infty             E[X^n]     \frac{t^n}{n!} - \sum_{n=1}^\infty                  n E[X^{n-1}] \frac{t^n}{n!}\\
&= E[X^0] + \sum_{n=1}^\infty  \Big( E[X^n] - n E[X^{n-1}] \Big)                \frac{t^n}{n!}.
\end{split}\end{equation}

Therefore, by \eqref{29} and \eqref{30}, we obtain the following theorem.

\begin{thm}
Let $X$ be a gamma random variable with parameters $\alpha=u$ and $\lambda=1$, and let $U$ be a uniform  random variable on (0,1). Assume that $X$ and $U$ are independent. For $n \geq 1$, we have
\begin{equation*}\begin{split}
 E[X^n] - n E[X^{n-1}] = (-1)^n b_n.
\end{split}\end{equation*}
\end{thm}

Let $X_1,X_2, \cdots X_k$ be independent gamma random variables with parameters 1,1, and let $U_1,U_2 \cdots ,U_k$ be uniformly independent random variables on (0,1). Assume that $X_i$ and $U_j$ are independent for all $i,j$. Then we have
\begin{equation}\begin{split}\label{31}
\sum_{n\geq k} S_{1,\lambda }(n,k) \frac{t^n}{n!} &= \frac{1}{k!} \Big( \log(1+\lambda t)^{\frac{1}{\lambda }}\Big)^k = \frac{\lambda ^{-k}}{k!} \Big( \log(1+\lambda t) \Big)^k \\
&=\frac{1}{k!} t^k \left( \sum_{l=0}^\infty \frac{(-1)^l}{l+1} \lambda ^l t^l \right)^k \\
&= \frac{1}{k!} t^k \left( \sum_{l_1=0}^\infty \frac{(-1)^{l_1}}{l_1+1} \lambda ^{l_1} t^{l_1} \right) \cdots \left( \sum_{l_k=0}^\infty \frac{(-1)^{l_k}}{l_k+1} \lambda ^{l_k} t^{l_k} \right)\\
&= \frac{1}{k!} t^k \sum_{n=0}^\infty \sum_{l_1+\cdots +l_k=n} (-1)^{l_1+\cdots+l_k} t^{l_1+\cdots+l_k}\\
&\quad \times \lambda ^{l_1+\cdots+l_k}  E[U_1^{l_1}] \cdots E[U_k^{l_k}]\\
&= \frac{1}{k!} t^k \sum_{n=0}^\infty \frac{(-\lambda )^n}{n!}t^n  \sum_{l_1+\cdots +l_k=n} {n \choose l_1,\cdots,l_k}\\
&\quad \times E[U_1^{l_1}\cdots U_k^{l_k}]l_1! \cdots l_k!\\
&=\frac{1}{k!} t^k \sum_{n=0}^\infty \frac{(-\lambda )^n}{n!}t^n  \sum_{l_1+\cdots +l_k=n} {n \choose l_1,\cdots,l_k}\\
&\quad \times E[U_1^{l_1}\cdots U_k^{l_k}] E[X_1^{l_1} \cdots E[X_k^{l_k}]\\
&=\frac{1}{k!} t^k \sum_{n=0}^\infty \frac{(-\lambda )^n}{n!}t^n \\
&\quad \times E\left[ \sum_{l_1+\cdots+l_k=n} {n \choose l_1,\cdots,l_k} (U_1X_1)^{l_1} \cdots (U_kX_k)^{l_k} \right]\\
&=\frac{t^k}{k!} \sum_{n=0}^\infty (-\lambda)^n E [ (U_1X_1 + \cdots U_kX_k)^n ] \frac{t^n}{n!}\\
&= \sum_{n=k}^\infty (-\lambda)^{n-k} E[(U_1X_1+\cdots+U_kX_k)^{n-k}] \frac{n!}{(n-k)! k!} \frac{t^n}{n!}\\
&=\sum_{n=k}^\infty (-\lambda)^{n-k} {n \choose k}  E[(U_1X_1+\cdots+U_kX_k)^{n-k}] \frac{t^n}{n!}.
\end{split}\end{equation}

Comparing the coefficients on both sides of \eqref{31}, we obtain the following theorem.

\begin{thm}
Let $X_1,X_2,\cdots ,X_k$ be independent gamma random variables with parameters 1,1, and let $U_1,\cdots,U_k$ be uniformly independent random variables on (0,1). Assume that $X_i$ and $U_j$ are independent for all $i$ and $j$. For $n,k \geq 0$ and $n \geq k$, we have
\begin{equation*}\begin{split}
S_{1,\lambda }(n,k) = (-\lambda)^{n-k}{n \choose k}  E[(U_1X_1+\cdots+U_1X_k)^{n-k}] .
\end{split}\end{equation*}
\end{thm}

Let $X_1,X_2,\cdots ,X_k$ be independent gamma random variables with parameters $\alpha=u$, $\lambda=1$ and $U$ be a uniform random variable on (0,1). Assume that $X_i$ and $U$ are independent for all $i$. From \eqref{28}, we have
\begin{equation}\begin{split}\label{32}
E[e^{(X_1+X_2+\cdots+X_k)t}] &= \left( \frac{-t}{(1-t) \log(1-t)} \right)^k = \left( \frac{-t}{\log(1-t)} \right)^k (1-t)^{-k} \\
&= \sum_{n=0}^\infty B_n^{(n-k+1)}(-k+1)(-1)^n \frac{t^n}{n!}.
\end{split}\end{equation}

On the other hand,
\begin{equation}\begin{split}\label{33}
E[e^{(X_1+X_2+\cdots+X_k)t}] &= \sum_{n=0}^\infty E[(X_1+\cdots+X_k)^n]                 \frac{t^n}{n!}.
\end{split}\end{equation}

Therefore, by \eqref{32} and \eqref{33}, we obtain the following theorem.

\begin{thm}
Let $X_1,X_2,\cdots ,X_k$ be independent gamma random variables with parameters $\alpha=u$, $\lambda=1$, and let $U$ be a uniform random variable on (0,1). Assume that $X_i$ and $U$ are independent for all $i$. For $n \geq 0$, we have
\begin{equation*}\begin{split}
 E[(X_1+\cdots+X_k)^n] = (-1)^n B_n^{(n-k+1)}(-k+1).
\end{split}\end{equation*}
\end{thm}

Assume that $U_1, U_2, \cdots ,U_k$ are uniformly independent random variables on (0,1). Then we observe that
\begin{equation}\begin{split}\label{34}
E[(1+\lambda t)^{\frac{U_1}{\lambda }}] &= \int_0^1 (1+\lambda t)^{\frac{u_1}{\lambda }} f(u_1) du_1 = \int_0^1 e^{\frac{u_1}{\lambda }\log(1+\lambda t)} du_1\\
&= \frac{1}{\frac{1}{\lambda }\log(1+\lambda t)} \big( (1+\lambda t)^{\frac{1}{\lambda }}-1\big).
\end{split}\end{equation}
From \eqref{34}, we note that
\begin{equation}\begin{split}\label{35}
E[(1+\lambda t)^{\frac{U_1+\cdots+U_k}{\lambda }}] &= E[(1+\lambda t)^{\frac{U_1}{\lambda }}]\times  \cdots \times E[(1+\lambda t)^{\frac{U_k}{\lambda }}] \\
&= \frac{k!}{\Big( \frac{1}{\lambda }\log(1+\lambda t)\Big)^k} \frac{1}{k!} \Big( (1+\lambda t)^{\frac{1}{\lambda }}-1 \Big)^k\\
&= \frac{k!}{\Big( \frac{1}{\lambda }\log(1+\lambda t)\Big)^k} \sum_{ n=k}^\infty S_{2,\lambda }(n,k) \frac{t^n}{n!}.
\end{split}\end{equation}
Thus, by \eqref{35}, we get
\begin{equation}\begin{split}\label{36}
\sum_{n=k}^\infty S_{2,\lambda }(n,k) \frac{t^n}{n!}& = \frac{\Big( \frac{1}{\lambda }\log(1+\lambda t)\Big)^k}{k!} \sum_{m=0}^\infty \left( \frac{\log(1+\lambda t)}{\lambda } \right)^m \frac{E[(U_1+\cdots+U_k)^m]}{m!} \\
&=\sum_{m=k}^\infty	\frac{m!}{(m-k)!k!} E[(U_1+\cdots+U_k)^{m-k}] \frac{1}{m!} \left( \frac{\log(1+\lambda t)}{\lambda } \right)^m\\
&= \sum_{m=k}^\infty E[(U_1+\cdots+U_k)^{m-k}] {m \choose k} \sum_{n=m}^\infty S_{1,\lambda }(n,m) \frac{t^n}{n!}\\
&= \sum_{n=k}^\infty  \left( \sum_{m=k}^n {m \choose k} S_{1,\lambda }(n,m) E[(U_1+\cdots+U_k)^{m-k}]             \right)   \frac{t^n}{n!}
\end{split}\end{equation}

Comparing the coefficients on both sides of \eqref{36}, we obtain the following theorem.

\begin{thm}
Let $U_1, U_2, \cdots, U_k$ be uniformly independent random variables on (0,1). For $n \in \mathbb{N}$ with $n \geq k$, we have
\begin{equation*}\begin{split}
 S_{2,\lambda }(n,k) = \sum_{m=k}^n {m \choose k} S_{1,\lambda }(n,m) E[(U_1+\cdots+U_k)^{m-k}]           .
\end{split}\end{equation*}
\end{thm}

Now, we observe that
\begin{equation}\begin{split}\label{37}
&E[(1+\lambda t)^{\frac{1}{\lambda }(U_1+U_2+\cdots+U_{k-1}+1)}] = (1+\lambda t)^{\frac{1}{\lambda }} E[(1+\lambda t)^{\frac{U_1+\cdots+U_{k-1}}{\lambda }}]\\
&=(1+\lambda t)^{\frac{1}{\lambda }} \Big(\frac{1}{\lambda }\log(1+\lambda t)\Big)^{-(k-1)} \Big((1+\lambda t)^{\frac{1}{\lambda }}-1\Big)^{k-1}.
\end{split}\end{equation}

Thus, by \eqref{37}, we get
\begin{equation}\begin{split}\label{38}
&\Big( \frac{1}{\lambda } \log(1+\lambda t)\Big)^{k-1} E[(1+\lambda t)^{\frac{1}{\lambda }(U_1+U_2+\cdots+U_{k-1}+1)}]\\
&=\Big( (1+\lambda t)^{\frac{1}{\lambda }}-1 \Big)^{k-1} \Big( (1+\lambda t)^{\frac{1}{\lambda }}-1+1) \\
&= \frac{k!}{k!}\Big( (1+\lambda t)^{\frac{1}{\lambda }}-1 \Big)^k + \frac{(k-1)!}{(k-1)!} \Big( (1+\lambda t)^{\frac{1}{\lambda }}-1 \Big)^{k-1}\\
&= k! \sum_{n=k}^\infty S_{2,\lambda }(n,k) \frac{t^n}{n!} + (k-1)! \sum_{n=k-1}^\infty S_{2,\lambda }(n,k-1) \frac{t^n}{n!}\\
&= (k-1)! \sum_{n=k-1}^\infty \Big( k S_{2,\lambda }(n,k) + S_{2,\lambda }(n,k-1) - n\lambda S_{2,\lambda }(n,k) + n\lambda S_{2,\lambda }(n,k) \Big) \frac{t^n}{n!}\\
&= (k-1)! \sum_{n=k-1}^\infty \Big( S_{2,\lambda }(n+1,k) + n\lambda  S_{2,\lambda }(n,k) \Big) \frac{t^n}{n!}.
\end{split}\end{equation}

On the other hand,
\begin{equation}\begin{split}\label{39}
&\Big( \frac{1}{\lambda } \log(1+\lambda t)\Big)^{k-1} E[(1+\lambda t)^{\frac{1}{\lambda }(U_1+U_2+\cdots+U_{k-1}+1)}]\\
&= \sum_{m=0}^\infty E[(U_1+U_2+\cdots+U_k+1)^m ] \frac{1}{m!} \Big( \frac{1}{\lambda } \log( 1+\lambda t)\Big)^{m+k-1}\\
&= \sum_{m=k-1}^\infty  E[(U_1+U_2+\cdots+U_k+1)^{m-k+1} ]\frac{m!}{(m-k+1)!} \sum_{n=m}^\infty S_{1,\lambda}(n,m) \frac{t^n}{n!}\\
&= \sum_{n=k-1}^\infty \left( \sum_{m=k-1}^n  E[(U_1+\cdots+U_k+1)^{m-k+1} ] S_{1,\lambda }(n,m) {m \choose k-1} (k-1)! \right) \frac{t^n}{n!}.
\end{split}\end{equation}

From \eqref{38} and \eqref{39}, we have
\begin{equation}\begin{split}\label{40}
 &S_{2,\lambda }(n+1,k) + n\lambda  S_{2,\lambda }(n,k) \\
 &= \sum_{m=k-1}^n  {m \choose k-1}S_{1,\lambda }(n,m) E[(U_1+\cdots+U_k+1)^{m-k+1} ] .
\end{split}\end{equation}
By replacing $n$ by $n-1$, we get
\begin{equation}\begin{split}\label{41}
&S_{2,\lambda }(n,k)+(n-1)\lambda S_{2,\lambda }(n-1,k)\\
&= \sum_{m=k-1}^{n-1} {m \choose k-1} S_{1,\lambda }(n-1,m) E[(U_1+\cdots +U_{k-1}+1)^{m-k+1}].
\end{split}\end{equation}

From Theorem 7 and \eqref{41}, we have
\begin{equation}\begin{split}\label{42}
&S_{2,\lambda }(n,k)+(n-1)\lambda S_{2,\lambda }(n-1,k)\\
&= \sum_{m=k}^n S_{1,\lambda }(n,m) {m \choose k} E[(U_1+\cdots+U_k)^{m-k}] \\
&\quad + (n-1)\lambda  \sum_{m=k}^{n-1} S_{1,\lambda }(n-1,m) {m \choose k} E[(U_1+\cdots+U_k)^{m-k}]\\
&={n \choose k} E[(U_1+\cdots+U_k)^{n-k}] S_{1,\lambda }(n,n) \\
&\quad + \sum_{m=k}^{n-1} E[(U_1+\cdots+U_k)^{m-k}] {m \choose k} S_{1,\lambda }(n-1,m-1).
\end{split}\end{equation}

On the other hand,
\begin{equation}\begin{split}\label{43}
&\sum_{m=k-1}^{n-1} {m \choose k-1} S_{1,\lambda }(n-1,m) E[(U_1+\cdots +U_{k-1}+1)^{m-k+1}]\\
&= {n-1 \choose k-1} S_{1,\lambda }(n-1,n-1) E[(U_1+\cdots+U_{k-1}+1)^{n-k}\\
&\quad + \sum_{m=k-1}^{n-2}{m \choose k-1} S_{1,\lambda }(n-1,m) E[(U_1+\cdots +U_{k-1}+1)^{m-k+1}].
\end{split}\end{equation}

From \eqref{41}, \eqref{42} and \eqref{43}, we have
\begin{equation}\begin{split}\label{44}
&{n \choose k} E[(U_1+\cdots+U_k)^{n-k}] S_{1,\lambda }(n,n) \\
&\quad- {n-1 \choose k-1} S_{1,\lambda }(n-1,n-1) E[(U_1+\cdots+U_{k-1}+1)^{n-k}\\
&=-\sum_{m=k-1}^{n-1} E[(U_1+\cdots+U_k)^{m-k}] {m \choose k} S_{1,\lambda }(n-1,m-1)\\
&\quad + \sum_{m=k-1}^{n-2}{m \choose k-1} S_{1,\lambda }(n-1,m) E[(U_1+\cdots +U_{k-1}+1)^{m-k+1}].
\end{split}\end{equation}

Thus, by \eqref{44}, we get
\begin{equation*}\begin{split}
&\frac{n}{k} E[(U_1+\cdots+U_k)^{n-k}] S_{1,\lambda }(n,n) - E[(U_1+\cdots+U_{k-1}+1)^{n-k}] S_{1,\lambda }(n-1,n-1)\\
&= \frac{1}{{n-1 \choose k-1}} \sum_{m=k-1}^{n-2} {m \choose k-1} S_{1,\lambda }(n-1,m) E[(U_1+\cdots+U_{k-1}+1)^{m-k+1}]\\
&\quad - \frac{1}{{n-1 \choose k-1}} \sum_{m=k-1}^{n-1} E[(U_1+\cdots+U_k)^{m-k} ] {m \choose k} S_{1,\lambda }(n-1,m-1).
\end{split}\end{equation*}

As $\lambda  \rightarrow 0$, we have
\begin{equation*}\begin{split}
\frac{n}{k} E[(U_1+\cdots+U_k)^{n-k}] = E[(U_1+\cdots+U_{k-1}+1)^{n-k}].
\end{split}\end{equation*}

A derangement is a permutation with no fixed points. For example, (2,3,1) and (3,1,2) are derangements of (1,2,3), but (3,2,1) is not because 2 is a fixed point. The number of derangements of an $n$-element set is called the $n$-th derangement number and denoted by $d_n$. This number satisfies the following recurrences:
\begin{equation}\begin{split}\label{45}
d_n = n \cdot d_{n-1} + (-1)^n,\,\,(n \geq 0),\quad (\textnormal{see} \,\, [3,4,7]).
\end{split}\end{equation}

By \eqref{45}, we get
\begin{equation}\begin{split}\label{46}
d_n= n! \sum_{k=0}^n \frac{(-1)^k}{k!},\,\,(n \geq 0),\quad (\textnormal{see} \,\, [7]).
\end{split}\end{equation}

From \eqref{46}, we can derive the following generating function.
\begin{equation}\begin{split}\label{47}
\frac{1}{1-t} e^{-t} &= \left( \sum_{k=0}^\infty \frac{(-1)^k}{k!} t^k \right)
\left( \sum_{m=0}^\infty t^m \right) = \sum_{n=0}^\infty \left( n! \sum_{k=0}^n \frac{(-1)^k}{k!} \right) \frac{t^n}{n!}\\
&=\sum_{n=0}^\infty            d_n      \frac{t^n}{n!},\quad (\textnormal{see} \,\, [3,4,7]).
\end{split}\end{equation}

Recently, the derangement polynomials are defined by the generating function
\begin{equation*}\begin{split}
\frac{1}{1-xt}e^{-t} =  \sum_{n=0}^\infty d_n(x) \frac{t^n}{n!},\quad (\textnormal{see} \,\, [7]).
\end{split}\end{equation*}
When $x=1$, $d_n(1) = d_n$, $(n \geq 0)$.

Let $X$ be gamma random variable with parameters 1,1. Then we have
\begin{equation}\begin{split}\label{48}
E[e^{Xt}] = \int_0^\infty e^{xt} e^{-t} dx = \frac{1}{1-t} = \frac{1}{1-t}e^{-t} e^t.
\end{split}\end{equation}

By \eqref{48}, we get
\begin{equation}\begin{split}\label{49}
E[e^{(X-1)t}] = \frac{1}{1-t} e^{-t} = \sum_{n=0}^\infty     d_n             \frac{t^n}{n!}.
\end{split}\end{equation}

From \eqref{49}, we have
\begin{equation*}\begin{split}
E[(X-1)^n ] = d_n, \,\,(n \geq 0).
\end{split}\end{equation*}

For $k \in \mathbb{N}$, we have
\begin{equation}\begin{split}\label{50}
E[e^{kXt}] = \int_0^\infty e^{kxt} e^{-x} dt = \left( \frac{1}{1-kt}e^{-t} \right) e^t.
\end{split}\end{equation}

Thus, by \eqref{50}, we get
\begin{equation}\begin{split}\label{51}
E[e^{(kX-1)t}] = \frac{1}{1-kt}e^{-t} =\sum_{n=0}^\infty            d_n(k)      \frac{t^n}{n!}.
\end{split}\end{equation}

By \eqref{51}, we get
\begin{equation}\begin{split}\label{52}
E[(kX-1)^n] = d_n(k),\,\,(n \geq 0).
\end{split}\end{equation}

Let $X_1,X_2,\cdots ,X_k$ be independent gamma random variables with parameters 1,1 . Then we have
\begin{equation}\begin{split}\label{53}
E[e^{(X_1+2X_2+\cdots+kx_k-k)t}] &= \left( \frac{1}{1-t} \right) \times  \left( \frac{1}{1-2t} \right) \times  \cdots  \left( \frac{1}{1-kt} \right) \times  e^{-kt}\\
&= \left( \frac{1}{1-t}e^{-t} \right) \times  \left( \frac{1}{1-2t}e^{-t} \right) \times  \cdots  \left( \frac{1}{1-kt} e^{-t}\right)\\
&= \left( \sum_{l_1=0}^\infty d_{l_1} \frac{t^{l_1}}{l_1!} \right) \times \left( \sum_{l_2=0}^\infty d_{l_2}(2) \frac{t^{l_2}}{l_2!} \right) \times \cdots \times \left( \sum_{l_k=0}^\infty d_{l_k}(k) \frac{t^{l_k}}{l_k!} \right)\\
&= \sum_{n=0}^\infty     \left( \sum_{l_1+\cdots+l_k=n} {n \choose l_1,\cdots,l_k} d_{l_1} d_{l_2}(2) \cdots d_{l_k}(k)          \right)\frac{t^n}{n!}.
\end{split}\end{equation}

On the other hand,
\begin{equation}\begin{split}\label{54}
E[e^{(X_1+2X_2+\cdots+kx_k-k)t}] = \sum_{n=0}^\infty E[(X_1+2X_2+\cdots+kX_k-k)^n]                  \frac{t^n}{n!}
\end{split}\end{equation}

Therefore, by \eqref{53} and \eqref{54}, we obtain the following theorem.
\begin{thm}
Let $X_1,X_2,\cdots, X_k$ be independent gamma random variables with parameters 1,1 . For $n \geq 0$, $k \in \mathbb{N}$, we have
\begin{equation*}\begin{split}
E[(X_1+2X_2+\cdots+kX_k-k)^n]=\sum_{l_1+\cdots+l_k=n} {n \choose l_1,\cdots,l_k} d_{l_1} d_{l_2}(2) \cdots d_{l_k}(k)  .
\end{split}\end{equation*}
\end{thm}

Now, we observe that
\begin{equation}\begin{split}\label{55}
\frac{1}{(1-t)(1-2t)\cdots(1-kt)} = \frac{1}{k!} \sum_{l=0}^k {k \choose l} (-1)^{k-l} l^k \frac{1}{1-lt}.
\end{split}\end{equation}

From \eqref{55}, we have
\begin{equation}\begin{split}\label{56}
&\frac{1}{k!} \sum_{l=0}^k {k \choose l} (-1)^{k-l} l^k \frac{1}{1-lt} = \sum_{m=0}^\infty \frac{1}{k!} \sum_{l=0}^k {k \choose l} (-1)^{k-l} l^{m+k} t^m\\
&= \sum_{m=0}^\infty \left( \frac{1}{k!} \Delta^k 0^{m+k} \right) t^m = \sum_{m=0}^\infty S_2(m+k,k) t^m,
\end{split}\end{equation}

where $\Delta f(x) = f(x+1)-f(x)$. By \eqref{56}, we easily get
\begin{equation}\begin{split}\label{57}
	&\frac{1}{(1-t)(1-2t)\cdots(1-kt)}  e^{-kt} = \left( \sum_{m=0}^\infty S_2(m+k,k) t^m \right) \left( \sum_{j=0}^\infty \frac{(-k)^j}{j!}t^j \right) \\
	&=\sum_{n=0}^\infty \left( \sum_{m=0}^n S_2(m+k,k) k^{n-m} (-1)^{n-m} \frac{n!}{(n-m)!} \right)                  \frac{t^n}{n!}\\
	&= \sum_{n=0}^\infty \left( \sum_{m=0}^n m! {n \choose m} (-1)^{n-m} S_2(m+k,k) k^{n-m} \right)                 \frac{t^n}{n!}.
	\end{split}\end{equation}	

Therefore, by \eqref{53} and \eqref{57}, we obtain the following theorem.

\begin{thm}
Let $X_1,X_2,\cdots ,X_k$ be independent gamma random variables with parameters 1,1 . For $n \geq 0$, $k \in \mathbb{N}$, we have
\begin{equation*}\begin{split}
E[(X_1+2X_2+\cdots+kX_k-k)^n ] =  \sum_{m=0}^n S_2(m+k,k) m! {n \choose m} (-1)^{n-m}  k^{n-m}.
\end{split}\end{equation*}
\end{thm}

{\noindent \bf{Remark.}} From Theorem 8 and Theorem 9, we have
\begin{equation*}\begin{split}
&\sum_{l_1+\cdots+l_k=n} {n \choose l_1,\cdots,l_k} d_{l_1} d_{l_2}(2) \cdots d_{l_k}(k) =\sum_{m=0}^n S_2(m+k,k) m! {n \choose m} (-1)^{n-m}  k^{n-m}.
\end{split}\end{equation*}

\end{document}